\documentclass[12pt]{article}
\usepackage[english]{babel}
\usepackage{amssymb, amsmath, amsfonts}
\usepackage{mathrsfs}
\newcommand{\p}{{\sf P}}
\newcommand{\e}{{\sf E}}
\newtheorem{Df}{Definition}
\newtheorem{Th}{Theorem}
\newtheorem{Lm}{Lemma}
\newtheorem{Rm}{Remark}

\begin{document}
\begin{center}
{\bf ON THE NEWMAN CONJECTURE}
\end{center}
\vskip0,5cm
\begin{center}
Alexander BULINSKI\footnote{The work is partially supported by RFBR grant 10-01-00397.}$^,$\footnote{Lomonosov
Moscow State University and University Paris-6 -- Pierre and Marie Curie.}
\end{center}
\vskip1cm

\begin{abstract}
We consider a random field, defined on an
integer-valued $d$-dimensional lattice $\mathbb{Z}^d$, with
covariance function satisfying a condition more general than
summability. Such condition appeared in the well-known Newman's
conjecture concerning the central limit theorem (CLT) for stationary
associated random fields. As was demonstrated by Herrndorf and
Shashkin, the conjecture fails already for $d=1$. In the present
paper, we show the validity of modified conjecture leaving intact
the mentioned condition on covariance function. Thus we establish,
for any integer $d\geq1$, a criterion of the  CLT validity  for the
wider class of positively associated stationary fields. The uniform
integrability for the squares of normalized partial sums, taken over
growing parallelepipeds or cubes in $\mathbb{Z}^d$, plays the key
role in deriving their asymptotic normality. So our result extends
the Lewis theorem proved for sequences of random variables. A
representation of variances of partial sums of a field using the
slowly varying functions in several arguments is employed in
essential way.
\vskip0,5cm
{\it Keywords and phrases}: stationary random fields, positive association, central limit theorem,
uniform integrability, slowly varying functions, the Newman conjecture.
\vskip0,5cm
{\it $2010$ AMS classification}: 60F05, 60G60.
\end{abstract}

\section{Introduction}

The study of asymptotical behavior of the
(normalized) sums of random variables is the vast research domain of
Probability Theory having various applications. The limit theorems
established for independent summands form here the classical core. In
this regard one can refer to the monographs \cite{GK}, \cite{Z},
\cite{IL}, \cite{P}; see also references therein.

Stochastic models described by means of families of dependent random variables
arose at the beginning of the last century.
Thus the Gaussian and Markov processes, martingales, solutions of the stochastic differential equations, mixing processes appeared
as well as other important classes (see, e.g., \cite{BS}, \cite{K}).  Moreover,
much attention has been paid to studying of random fields.

Since the 1960s due to the problems of mathematical statistics, reliability theory, percolation and
statistical physics there arose
the stochastic models based on the families of variables possessing various forms of positive
or negative dependence (see, e.g., \cite{BSha}).
The key role in these models belongs to the notion of
association (in statistical physics the well-known FKG-inequalities imply the association).
We will use the following concept
extending that introduced in \cite{EPW}.

\begin{Df}\label{PA}$(${\rm \cite{N}}$)$
A real-valued random field $X=\{X_t,t\in T\}$ is called positively associated\footnote{or weakly associated} $($one writes $X\in {\sf PA})$
if, for any finite disjoint sets $I=\{s_1,\ldots,s_m\}\subset T$, $J=\{t_1,\ldots,t_n\}\subset T$ and all bounded coordinate-wise nondecreasing Lipschitz functions $f:\mathbb{R}^m\to \mathbb{R}$, $g:\mathbb{R}^n\to \mathbb{R}$, one has
\begin{equation}\label{e1}
{\sf cov}(f(X_{s_1},\ldots,X_{s_m}),g(X_{t_1},\ldots,X_{t_n})) \geq 0.
\end{equation}
\end{Df}

Recall that a random field $X$ is called associated (\cite{EPW}), if the definition above is satisfied without the hypothesis $I\cap J=\varnothing$. Obviously association implies positive association. Note that any family of (real-valued) independent random variables
is automatically associated. Many other important examples can be found in \cite{BSha}.

For a random field $X=\{X_t,t\in T\}$ and a finite set $U\subset T$ introduce
$$
S(U)=\sum_{t\in U} X_t.
$$
Further on we will consider random fields defined on a lattice $T=\mathbb{Z}^d$ and a probability space $(\Omega,\mathcal{F},\p)$.
In the seminal paper by Newman \cite{N} the central limit theorem (CLT) was established for associated (strictly) stationary random
field $X=\{X_t,t\in \mathbb{Z}^d\}$ under {\it finite susceptibility condition} that is when the covariance function is summable:
\begin{equation}\label{newm}
\sigma^2:=\sum_{j\in\mathbb{Z}^d} {\sf cov}(X_0,X_j) < \infty.
\end{equation}
Namely, these simple assumptions imply for a field $X$ the following relation
\begin{equation}\label{clt0}
\frac{S_n - {\sf E}S_n}{\sqrt{\langle n \rangle}} \stackrel{law}\longrightarrow Z \sim \mathcal{N}(0,\sigma^2)\;\;\mbox{as}\;\;n\to \infty,\;n=(n_1,\ldots,n_d)\in\mathbb{N}^d,
\end{equation}
here $S_n = S([0,n]\cap \mathbb{Z}^d)$, $[0,n]=
[0,n_1]\times \ldots \times [0,n_d]$, $\langle n \rangle= n_1\ldots n_d$, $\mathcal{N}(0,\sigma^2)$
is a Gaussian law with parameters $0$ and $\sigma^2$,
$\stackrel{law}\longrightarrow$ stands for weak convergence of distributions.

The goal of this work is to provide the criteria of the CLT validity for positively associated stationary random fields with finite second moment (and in general without condition \eqref{newm}).

\section{Main results}

At first it is reasonable to recall several definitions.
\begin{Df}
A function $L:\mathbb{R}^d_+ \to \mathbb{R}\setminus \{0\}$ is called {\it slowly varying} $($at infinity$)$
if, for any vector  $a=(a_1,\ldots,a_d)^{\top}$ with positive coordinates,
\begin{equation}\label{svf}
\frac{L(a_1x_1,\ldots,a_d x_d)}{L(x_1,\ldots,x_d)} \to 1\;\;\mbox{as}\;\;x=(x_{1},\ldots,x_{d})^{\top}\to \infty,
\end{equation}
i.e. $x_1\to \infty,\ldots,x_d\to \infty$. For such functions we write $L\in \mathcal{L}(\mathbb{R}^d_+)$.
\end{Df}

We operate with column vectors and use the symbol $\top$ for transposition. A function  $L:\mathbb{N}^d \to \mathbb{R}\setminus \{0\}$
is called {\it slowly varying} (at infinity)
if, for any vector  $a=(a_1,\ldots,a_d)^{\top}\in \mathbb{N}^d$, relation \eqref{svf} holds with additional assumption that
$x\in  \mathbb{N}^d$. Then we write $L \in \mathcal{L}(\mathbb{N}^d)$.

For example the function $\prod_{k=1}^d \log(x_k\vee 1)$ where $x\in \mathbb{R}^d_+$ belongs to $ \mathcal{L}(\mathbb{R}^d_+)$.

\begin{Rm}\label{rm1} It is well-known that not every function belonging to $\mathcal{L}(\mathbb{N}^d)$ admits extension
to a function from the class $ \mathcal{L}(\mathbb{R}^d_+)$ even for $d=1$. However, it is not difficult to
verify that
if a coordinate-wise nondecreasing function $L\in \mathcal{L}(\mathbb{N}^d)$, then $H(x):=L([\widetilde{x}])$ belongs to
 $\mathcal{L}(\mathbb{R}^d_+)$. Here $\widetilde{x}=(x_1\vee 1,\ldots,x_d\vee 1)^{\top}$ for $x\in\mathbb{R}^d$, and
$[x]=([x_1],\ldots,[x_d])^{\top}$,
i.e. one takes the integer part of each component of $x$.
\end{Rm}

During a long time there was no solution to the Newman conjecture on possible replacement of requirement
\eqref{newm} appearing in CLT by a milder condition. Namely, he considered the partial sums $S(U)$ taken over "integer cubes" $U$ and
believed that instead of \eqref{newm} it suffices to assume that for associated strictly stationary random field
$X=\{X_j,j\in \mathbb{Z}^d\}$ with $\e X_0^2 < \infty$ the function
\begin{equation}\label{fn}
{\sf K}(r)= \sum_{j\in \mathbb{Z}^d:\|j\| \leq r}{\sf cov}(X_0,X_j),\;\;\;r\in\mathbb{N},
\end{equation}
belongs to $\mathcal{L}(\mathbb{N})$ where $\|\cdot\|$ is the Euclidean norm in $\mathbb{R}^d$.

Unfortunately it turned out that this beautiful hypothesis is not true even for $d=1$.
The first counterexample was constructed by Herrndorf \cite{H}, and then Shashkin \cite{S} showed that
condition \eqref{newm} has in a sense the optimal character.

It is worth mentioning also that the Newman CLT was  generalized in \cite{BV} for partial sums
$S(U)$ taken over regularly growing subsets of  $\mathbb{Z}^d$. Further extensions are discussed in
Chapter 3 of \cite{BSha}.

\begin{Df}\label{d3} A family $X=\{X_j,j\in\mathbb{N}^d\}$ is called uniformly integrable if
$$
\lim_{c\to \infty} \sup_{j\in \mathbb{N}^d} {\sf E}|X_j|\mathbb{I}\{|X_j|\geq c\} = 0.
$$
\end{Df}

For a (wide sense) stationary random field $X=\{X_j,j\in\mathbb{Z}^d\}$ introduce the function
$$
{\sf K}_X(n)=\sum_{j\in\mathbb{Z}^d: -n\leq j\leq n} {\sf cov}\,(X_0,X_j),\;\;n\in\mathbb{N}^d.
$$
If $a=(a^{(1)},\ldots,a^{(d)})^{\top}$ and $b= (b^{(1)},\ldots,b^{(d)})^{\top}$ are vectors in
$\mathbb{R}^d$, the notation $a\leq b$ means that  $a^{(k)} \leq b^{(k)}$ for all $k=1,\ldots,d$.
We write $a< b$ whenever $a^{(k)} < b^{(k)}$ for any $k=1,\ldots,d$.

The following result extends the Lewis theorem  proved in \cite{L} for a sequence of random variables.

\begin{Th}\label{th1}
Let a strictly stationary random field $X=\{X_j, j\in\mathbb{Z}^d\} \in {\sf PA}$,
$0<\e X_0^2 < \infty$ and ${\sf K}_X(\cdot) \in \mathcal{L}(\mathbb{N}^d)$. Then $X$ satisfies CLT,
i.e.
\begin{equation}\label{clt}
\frac{S_n - \e S_n}{\sqrt{{\sf var}S_n}} \stackrel{law}\longrightarrow Z \sim \mathcal{N}(0,1)\;\;\mbox{as}\;\;n\to \infty,
\end{equation}
if and only if the family $\{(S_n - \e S_n)^2/(\langle n \rangle {\sf K}_X(n)), n\in \mathbb{N}^d\}$ is uniformly integrable.
\end{Th}

Consider now a sequence of growing "integer cubes"  $C_r=(0,r]^d\cap \mathbb{Z}^d$,
$r\in\mathbb{N}$.

\begin{Th}\label{th2} Let a strictly stationary random field $X=\{X_j, j\in\mathbb{Z}^d\} \in {\sf PA}$,
$0<\e X_0^2 < \infty$ and ${\sf K}(\cdot) \in \mathcal{L}(\mathbb{N})$. Then
\begin{equation}\label{cltcubes}
\frac{S(C_r) - \e S(C_r)}{\sqrt{{\sf var}S(C_r)}} \stackrel{law}\longrightarrow Z
\sim \mathcal{N}(0,1)\;\;\mbox{as}\;\;r\to \infty,
\end{equation}
if and only if the sequence  $((S(C_r) - \e S(C_r))^2/(r^d\, {\sf K}(r)))_{r\in \mathbb{N}}$ is uniformly
integrable.
\end{Th}

Theorem \ref{th2} shows what one has to assume additionally, for a class of positively associated strictly stationary
random fields, besides the condition that the function ${\sf K}(\cdot)$ is slowly varying to guarantee that the
Newman conjecture holds true for any dimension $d\in\mathbb{N}$. In \cite{N1} the author discussed his conjecture and
noted without proof that the "mild version" of that hypothesis takes place under the additional condition of
uniform integrability of the sequence appearing in Theorem 2 above. Therefore Theorems \ref{th1} and \ref{th2} show that
in fact we do not change the initial problem but clarify its essential feature. We do not deal here with a
renorm group approach
(do not consider the partition of $\mathbb{R}^d$
by the congruent cubes) but study the partial sums $S_n$ taken over any growing "integer blocks".

\section{Proofs of the main results}

We start with  simple auxiliary statements.

\begin{Lm}\label{lm1}
Let a function $L$ belonging to $\mathcal{L}(\mathbb{N}^d)$ be coordinate-wise nondecreasing. Then there exist
non-random vectors $q_n=(q_n^{(1)},\ldots,q_n^{(d)})^{\top} \in\mathbb{N}^d$, where $n
=(n_1,\ldots,n_d)^{\top}\in\mathbb{N}^d$, such that
\begin{equation}\label{eq1}
q_n^{(k)} \leq n_k,\;\frac{q_n^{(k)}}{n_k} \to 0\;\;\mbox{for}\;\;k=1,\ldots,d,\;\;q_n\to \infty\;\;\;\mbox{and}\;\;\;
\frac{L(n)}{L(q_n)}\to 1\;\mbox{as}\;n\to \infty.
\end{equation}
\end{Lm}

{\bf Proof.} According to Remark \ref{rm1}
we can assume without loss of generality that $L$ is extended to a function belonging to
the class $ \mathcal{L}(\mathbb{R}^d_+)$. For any  $R=(R^{(1)},\ldots,R^{(d)})^{\top} \in \mathbb{N}^d$
we can choose $N_0(R)\in\mathbb{N}^d$ in such a way that
$$
\frac{L(n_1,\ldots,n_d)}{L\Bigl(\frac{n_1}{R^{(1)}},\ldots,\frac{n_d}{R^{(d)}}\Bigr)}-1\leq \frac{1}{\langle R \rangle}
$$
for all  $n\geq N_0(R)$. Now we take a sequence $(R(r))_{r\in\mathbb{N}}$ such that $R(r) \in \mathbb{N}^d$ and
$R(r)<R(r+ 1)$ for each $r\in \mathbb{N}$. Introduce
$M_0(1)=N_0(R(1))$ and $M_0(r+1)= (M_0(r)\vee N_0(R(r+1)))+{\bf 1}$
for $r\in \mathbb{N}$ where, as usual,  ${\bf 1}=(1,\ldots,1)^{\top} \in \mathbb{R}^d$ and
$$
(a^{(1)},\ldots,a^{(d)})\vee (b^{(1)},\ldots,b^{(d)})= (a^{(1)}\vee b^{(1)},\ldots,a^{(d)}\vee b^{(d)}).
$$
Then $M_0(r)<M_0(r+1)$ for $r\in\mathbb{N}$. For arbitrary $r\in\mathbb{N}$ and $n\geq M_0(r)$
$$
\frac{L(n_1,\ldots,n_d)}{L\Bigl(\frac{n_1}{R^{(1)}(r)},\ldots,\frac{n_d}{R^{(d)}(r)}\Bigr)}-1\leq \frac{1}{\langle R(r) \rangle}.
$$
Let us define non-random sequences  $(\varepsilon_j^{(k)})_{j\in \mathbb{N}}$ where $k=1,\ldots,d$, putting
$\varepsilon_{j}^{(k)} = 1/R^{(k)}(r)$ for $M^{(k)}_0(r)\leq j < M^{(k)}_0(r+1)$.

For any $\varepsilon >0$ take $r_0 \in\mathbb{N}$ in such a way that $1/\langle R(r_0)\rangle<\varepsilon$.
Further on, for $n$ such that $M_0(r)\leq n<M_0(r+1)$ where $r\geq r_0$, one has
$$
1\leq \frac{L(n_1,\ldots,n_d)}{L(n_1\varepsilon_{n_1}^{(1)},\ldots,n_d\varepsilon_{n_d}^{(d)})}
=\frac{L(n_1,\ldots,n_d)}{L\Bigl(\frac{n_1}{R^{(1)}(r)},\ldots,\frac{n_d}{R^{(d)}(r)}\Bigr)}
$$
$$
\leq 1+\frac{1}{\langle R(r) \rangle} \leq 1+\frac{1}{\langle R(r_0) \rangle}\leq 1+\varepsilon.
$$
Then we can take
$
q_n=([n_1\varepsilon_{n_1}^{(1)}],\ldots,[n_d\varepsilon_{n_d}^{(d)}])\vee ([\log n_1],\ldots,[\log n_d])
\vee {\bf 1},
$
to ensure
 the validity of \eqref{eq1}. $\square$

\begin{Lm}\label{lm2}
Let $X=\{X_j, j\in\mathbb{Z}^d\}$ be a wide sense stationary random field with nonnegative covariance function.
Assume that ${\sf K}_X(\cdot)\in \mathcal{L}(\mathbb{N}^d)$. Then
\begin{equation}\label{vps}
{\sf var} S(U_n) \sim \langle n \rangle\,{\sf K}_X(n) \;\;\;\mbox{as}\;\;n\to \infty
\end{equation}
where $U_n=\{j\in\mathbb{Z}^d: {\bf 1}\leq j\leq n\}$, $n\in\mathbb{N}^d$. Conversely, if
${\sf var} S(U_n) \sim \langle n \rangle\, L(n)$
as $n\to \infty$, where $L\in\mathcal{L}(\mathbb{N}^d)$, then
$L(n)\sim {\sf K}_X(n)$ as $n\to \infty$.
\end{Lm}

{\bf Proof.} Let ${\sf K}_X(\cdot)\in \mathcal{L}(\mathbb{N}^d)$. Due to the (wide-sense) stationarity of $X$ one has
${\sf cov}(X_i,X_j)={\sf R}(i-j)$ for $i,j\in\mathbb{Z}^d$. Thus
$$
{\sf var} S(U_n) = \sum_{i,j\in U_n} {\sf cov}(X_i,X_j) = \sum_{i,j\in U_n}{\sf R}(i-j)
$$
$$
=\sum_{m\in \mathbb{Z}^d: -(n-{\bf 1})\leq m\leq n-{\bf 1}} (n_1 -|m_1|)\ldots (n_d -|m_d|){\sf R}(m)
$$
\begin{equation}\label{aux1a}
\leq \langle n \rangle \, \sum_{m\in \mathbb{Z}^d: -(n-{\bf 1})\leq m\leq n-{\bf 1}} {\sf R}(m)
\leq  \langle n \rangle \, {\sf K}_X(n),
\end{equation}
as the function ${\sf R}$ is nonnegative.

Take any $c\in (0,1)$ and $n \geq \frac{1}{1-c}{\bf 1}$
(i.e. $cn\leq n-{\bf 1}$, $n\in \mathbb{N}^d$). Using again nonnegativity of ${\sf R}$ we can write
$$
{\sf var} S(U_n)
=\sum_{m\in \mathbb{Z}^d: -(n-{\bf 1})\leq m\leq n-{\bf 1}} (n_1 -|m_1|)\ldots (n_d -|m_d|){\sf R}(m)
$$
$$
\geq (1-c)^d \langle n \rangle \,\sum_{m\in \mathbb{Z}^d: -cn \leq m \leq cn} {\sf R}(m) = (1-c)^d
\langle n \rangle \,{\sf K}_X([cn]).
$$
In view of Remark \ref{rm1} we come to the relation
$$
(1-c)^d \langle n \rangle \,{\sf K}_X([cn])\sim (1-c)^d \langle n \rangle \,{\sf K}_X(n),
\;\;n\to \infty, \;\;n\in\mathbb{N}^d.
$$
Consequently, ${\sf var}S(U_n) \sim \langle n \rangle \,{\sf K}_X(n)$ as $n\to \infty$, because
$c$ can be taken arbitrary close to zero.

Now suppose that ${\sf var} S(U_n)\sim \langle n \rangle \, L(n)$ as $n\to \infty$, where $L
\in \mathcal{L}(\mathbb{N}^d)$. Then for any $\varepsilon >0$ and all $n$ sufficiently large (i.e.
each component of $n$ is large enough), application of \eqref{aux1a} leads to the inequality
\begin{equation}\label{aux1}
{\sf K}_X(n) \geq \frac{{\sf var}S(U_n)}{ \langle n \rangle} \geq (1-\varepsilon)L(n).
\end{equation}
For a fixed $q\in\mathbb{N}$, $q>1$, and $n_r\in\mathbb{N}$, $m_r\in\mathbb{Z}$ such that
$|m_r|\leq n_r$ where $r=1,\ldots,d$, one has
$$
\frac{q}{q-1}\left(1-\frac{|m_r|}{n_rq}\right) \geq \frac{q}{q-1}\left(1-\frac{n_r}{n_rq}\right)=1.
$$
Therefore, taking into account condition ${\sf R}\geq 0$ we verify that
$$
{\sf K}_X(n) \leq \left(\frac{q}{q-1}\right)^d \sum_{m\in\mathbb{Z}^d: -n\leq m\leq n} {\sf R}(m) \prod_{r=1}^d \frac{(n_rq-|m_r|)}{n_rq}
$$
$$
\leq\left(\frac{q}{q-1}\right)^d\left(\prod_{r=1}^d n_rq\right)^{-1} \sum_{m\in\mathbb{Z}^d: -nq \leq m\leq nq} {\sf R}(m) \prod_{r=1}^d(n_rq-|m_r|)
$$
\begin{equation}\label{aux2}
=\left(\frac{q}{q-1}\right)^d\frac{{\sf var}S(U_{qn})}{\langle qn \rangle}\sim \left(\frac{q}{q-1}\right)^d L(qn),\;\;n\to \infty.
\end{equation}
As $q$ can be chosen arbitrary large, using \eqref{aux1} and \eqref{aux2} we conclude that the
desired statement holds.
$\square$

\vskip0,5cm
{\bf Proof of Theorem 1.} {\sf Necessity.} Suppose that \eqref{clt} is satisfied.
Then
$$
\frac{(S_n-\e S_n)^2}{{\sf var} S_n} \stackrel{law}\longrightarrow Z^2
\;\;\mbox{as}\;\;n\to \infty.
$$
Indeed, if the random variables $Y_n \stackrel{law}\longrightarrow  Y$, then for any bounded continuous function
$h:\mathbb{R}\to\mathbb{R}$ one has $h(Y_n)\stackrel{law}\longrightarrow h(Y)$ as $n\to \infty$.
Obviously,
$$
\frac{(S_n-\e S_n)^2}{{\sf var} S_n}\geq 0\;\;\mbox{and}\;\;\frac{\e(S_n-\e S_n)^2}{{\sf var} S_n}=1.
$$
Thus uniform integrability of the family $\{(S_n-\e S_n)^2/{\sf var} S_n, n\in\mathbb{N}^d\}$
follows from the analogue of Theorem 1.5.4 established in \cite{B} for a sequence  of random variables indexed by
points of $\mathbb{N}$. In view of Lemma \ref{lm2} we can claim that \eqref{vps} holds. Consequently, the family
$\{(S_n - \e S_n)^2/(\langle n \rangle {\sf K}_X(n)), n\in \mathbb{N}^d\}$ is also uniformly integrable.

{\sf Sufficiency.}
If the function ${\sf K}_X$ is bounded we see that \eqref{newm} is valid and Theorem 3.1.12 of \cite{BSha}
implies that \eqref{clt} is satisfied. Thus we will assume further that a function  ${\sf K}_X$ is unbounded.
Set
${\sf K}_X(t):={\sf K}_X([t]\vee {\bf 1})$
for $t=(t_1,\ldots,t_d)^{\top}\in \mathbb{R}^d_+$ where $[t]=([t_1],\ldots,[t_d])^{\top}$.
This extension of the initial function ${\sf K}_X$ belongs to $\mathcal{L}(\mathbb{R}^d_+)$ as
${\sf K}_X$ is coordinate-wise nondecreasing on $\mathbb{N}^d$ (a field $X\in {\sf PA}$, therefore its covariance function is nonnegative).
Further on we assume that the function ${\sf K}_X$ is extended on $\mathbb{R}^d_+$ as indicated above.

Let the vectors $q_n$, $n\in\mathbb{N}^d,$ be constructed according to Lemma \ref{lm1}.
It is not difficult to find a non-random family of vectors
$\{p_n, n\in\mathbb{N}^d\}$,
where $p_n$ takes values in $\mathbb{N}^d$, such that
\begin{equation}\label{eq2}
q_n^{(k)} \leq p_n^{(k)}\leq n_k,\;q_n^{(k)}/p_n^{(k)} \to 0\;\;\mbox{and}\;\;p_n^{(k)}/n_k \to 0\;\;
\mbox{for}\;k=1,\ldots,d\;\;\mbox{as}\;n\to \infty.
\end{equation}

Now we apply the Bernstein partitioning method. For $n,j\in\mathbb{N}^d$ and introduced $p_n$ and $q_n$ consider the blocks
$$
U_n^{(j)}=\{u\in \mathbb{N}^d: (j_k-1)(p_n^{(k)}+q_n^{(k)})<u_k \leq j_kp_n^{(k)}+(j_k -1)q_n^{(k)},\;k=1,\ldots,d\},
$$
where $u=(u_1,\ldots,u_d)$. Let $J_n=\{j\in\mathbb{N}^d: U_n^{(j)} \subset U_n\}$ and
$$
W_n =\bigcup_{j\in J_n} U_n^{(j)}, \;\;\;G_n= U_n\setminus W_n,\;\;n\in\mathbb{N}^d.
$$
In other words $W_n$ consists of "large blocks" (having the "size" $p_n^{(k)}$ along each of the $k-$th axis for $k=1,\ldots,d$),
separated by "corridors" belonging to the set $G_n$. Put  $v_n = \sqrt{\langle n \rangle {\sf K}_X(n)}$. Then, for
each $t\in\mathbb{R}$ and
$n\in\mathbb{N}^d$ we obtain
$$
\left|{\sf E}\exp\left\{\frac{it}{v_n}S_n\right\} - e^{-\frac{t^2}{2}}\right| \leq
\left|{\sf E}\exp\left\{\frac{it}{v_n}S_n\right\} -
{\sf E}\exp\left\{\frac{it}{v_n}\sum_{j\in J_n}S(U_n^{(j)})\right\}\right|
$$
$$
+\left|{\sf E}\exp\left\{\frac{it}{v_n}\sum_{j\in J_n}S(U_n^{(j)})\right\} - \prod_{j\in J_n}
{\sf E}\exp\left\{\frac{it}{v_n}S(U_n^{(j)})\right\}\right|
$$
$$
+\left|\prod_{j\in J_n}{\sf E}
\exp\left\{\frac{it}{v_n}S(U_n^{(j)}\right\}
- e^{-\frac{t^2}{2}}\right| =:\sum_{r=1}^3 Q_r,
$$
here $i^2=-1$, $Q_r=Q_r(n,t)$ and $S_n = S(U_n)$ as previously.
Taking into account that $|e^{ix}-e^{iy}|\leq |x-y|$ for all $x,y \in \mathbb{R}$,
and using the Lyapunov inequality we get
$$
Q_1 \leq \frac{|t|}{v_n} {\sf E}|S(G_n)| \leq  \frac{|t|}{v_n} ({\sf E}S(G_n)^2)^{1/2}.
$$
A random field $X\in {\sf PA}$, therefore ${\sf cov}\,(X_j,X_u) \geq 0$ for any $j,u\in \mathbb{N}^d$.
Thus in view of wide-sense stationarity of $X$ we come to the relations
$$
{\sf E} S(G_n)^2 \leq \sum_{j\in G_n} \sum_{u:-n\leq u-j \leq n} {\sf cov} (X_j,X_u) \leq
{\sf card}\, G_n \,{\sf K}_X(n)
$$
$$
\leq {\sf K}_X(n) \sum_{k=1}^d (m_n^{(k)}q_n^{(k)} + p_n^{(k)} + q_n^{(k)})
\prod_{1\leq l\leq d, l\neq k} n_l
$$
where ${\sf card}\, G$ stands for the cardinality of a set $G$, $m_n^{(k)} = [n_k/(p_n^{(k)}+ q_n^{(k)})]$, $k=1,\ldots,d$. Due to \eqref{eq1} and \eqref{eq2} we get the inequality
$$
\frac{{\sf E}S(G_n)^2}{\langle n \rangle\, {\sf K}_X(n)}
\leq \sum_{k=1}^d \frac{m_n^{(k)}q_n^{(k)} + p_n^{(k)} + q_n^{(k)}}{n_k}\to 0,\;\;\;n\to \infty.
$$
Consequently, $Q_1(n,t)\to 0$ for each $t\in\mathbb{R}$ as $n\to \infty$.

For any $n\in \mathbb{N}^d$ the family $\{S(U_n^{(j)}), j\in J_n\} \in {\sf PA}$ (see, e.g., Theorem 1.1.8
in \cite{BSha}). Enumerate elements of this family to obtain the collection of random variables
$\{Y_{n,s}, s=1,\ldots,M_n\}$ where $M_n= {\sf card}\,J_n$.
It is easily seen that
$$
\prod_{k=1}^d m_n^{(k)}
\leq M_n \leq \prod_{k=1}^d (m_n^{(k)}+1).
$$
Recall that for complex-valued random variables $Y$ and $V$
(absolute square integrable) the covariance ${\sf cov}(Y,V):= {\sf E}(Y-{\sf E}Y)
\overline{(V-{\sf E}V)}$, where the bar denotes the conjugation.
Due to Theorem 1.5.3 of \cite{BSha} one has
$$
Q_2 \leq \sum_{s=1}^{M_n-1}\left|\,{\sf cov}\!\left(\exp\left\{\frac{it}{v_n}Y_{n,s}\right\},\exp\left\{-\frac{it}{v_n}\sum_{l=s+1}^{M_n} Y_{n,l}\right\}\right)\right|
$$
$$
\leq \frac{4t^2}{v_n^2}\sum_{1\leq s,l \leq M_n,s\neq l}{\sf cov}(Y_{n,s},Y_{n,l})
\leq \frac{4t^2}{\langle n \rangle \,{\sf K}_X(n)}\sum_{j\in U_n}\sum_{u\in U_n,|u-j|>q_n}{\sf cov}(X_j,X_u)
$$
where $|u|=\max_{k=1,\ldots,d}|u_k|$.
Obviously, for $j\in U_n$
$$
\{u\in U_n, |u-j|>q_n\} \subset \{u\in \mathbb{N}^d: j-n \leq u \leq j+n\}\setminus \{u\in \mathbb{N}^d:|u-j|\leq q_n\}.
$$
Therefore, the inequality
$$
\sum_{j\in U_n} \sum_{u\in U_n, |u-j|>q_n} {\sf cov}(X_j,X_u) \leq \langle n \rangle ({\sf K}_X(n) - {\sf K}_X(q_n))
$$
and \eqref{eq2} imply that $Q_2(n,t)\to 0$  for each $t\in \mathbb{R}$ as $n\to \infty$.

For any $n\in\mathbb{N}$ introduce a vector $(Z_{n,1},\ldots,Z_{n,M_n})^{\top}$ having the independent components and such that
the law of $Z_{n,k}$ coincides with the law of $Y_{n,k}/v_n$, $k=1,\ldots,d$.   Due to Lemma \ref{lm2} for all
$s=1,\ldots,M_n$
\begin{equation}\label{aux3}
{\sf var} Z_{n,s} = {\sf var } Z_{n,1} \sim \langle p_n \rangle \, {\sf K}_X(p_n)/\langle n \rangle \,{\sf K}_X(n),\;\;\;n\to \infty.
\end{equation}
Thus
\begin{equation}\label{aux4}
\sum_{s=1}^{M_n} {\sf var} Z_{n,s} = M_n {\sf var} Z_{n,1} \to 1, \;\;\;n\to \infty,
\end{equation}
since
$$M_n\langle p_n \rangle \sim \prod_{k=1}^d [n_k/(p_n^{(k)} + q_n^{(k)})]p_n^{(k)} \sim \langle n \rangle $$
and ${\sf K}_X(p_n)/{\sf K}_X(n) \to 1$ as $n\to \infty$.
For arbitrary $\varepsilon >0$, taking into account the stationarity of $X$, we have
$$
\sum_{s=1}^{M_n} {\sf E}Z_{n,s}^2 \mathbb{I}\{|Z_{n,s}|>\varepsilon\} = \frac{M_n}{\langle n \rangle\, {\sf K}_X(n)} {\sf E} Y_{n,1}^2\mathbb{I}\{Y_{n,1}^2>\varepsilon^2 \langle n \rangle\, {\sf K}_X(n)\}
$$
$$
=\frac{M_n\langle p_n \rangle\, {\sf K}_X(p_n)}{\langle n \rangle\, {\sf K}_X(n)} {\sf E}\frac{S(U_n^{(1)})^2}{\langle p_n \rangle\, {\sf K}_X(p_n)}\mathbb{I}
\left\{
\frac{S(U_n^{(1)})^2}{\langle p_n \rangle\, {\sf K}_X(p_n)} > \varepsilon^2
\frac{\langle n \rangle\, {\sf K}_X(n)}
{\langle p_n \rangle\, {\sf K}_X(p_n)}
\right\} \to 0,\;\;n\to\infty,
$$
in view of  \eqref{aux3}, \eqref{aux4} and because
$$
\frac{\langle n \rangle {\sf K}_X(n)}{\langle p_n \rangle {\sf K}_X(p_n)}\to \infty\;\;\mbox{as}\;\;
n\to \infty.
$$
We also used uniform integrability of $\{S(U_n^{(1)})^2/(\langle p_n \rangle\, {\sf K}_X(p_n)), n\in \mathbb{N}^d\}$. Indeed, 
this is a subfamily of the uniformly integrable family
$\{S_n^2/(\langle n \rangle\, {\sf K}_X(n)), n\in \mathbb{N}^d\}$.
The Lindeberg theorem (see, e.g., \cite{K}, p. 69) implies that
$$
\sum_{s=1}^{M_n} Z_{n,s} \stackrel{law}\to Z\sim \mathcal{N}(0,1),\;\;\;n\to\infty.
$$
Therefore,
$$
\prod_{s=1}^{M_n}{\sf E} \exp\{itZ_{n,s}\} - \exp\left\{-\frac{t^2}{2}\right\} \to 0,\;\;n\to \infty.
$$
It remains to note that
$$
\prod_{j\in J_n} {\sf E}\exp\left\{\frac{it}{v_n}S(U^{(j)}_n)\right\} = \prod_{s=1}^{M_n}{\sf E}\exp\{itZ_{n,s}\}.
$$
Thus $Q_3(n,t)\to 0$  for each $t\in\mathbb{R}$ as $n\to\infty$.
The proof is complete. $\square$

\vskip0,5cm
{\bf Proof of Theorem 2.} For a wide-sense stationary random field  $X=\{X_j,j\in\mathbb{Z}^d\}$ introduce the function
$$
{\sf R}_X(r) = \sum_{j\in\mathbb{Z}^d: |j|\leq r}{\sf cov}(X_0,X_j), \;\;r\in \mathbb{N}.
$$
This function ${\sf R}_X(\cdot)$ is close in a sense to ${\sf K}(\cdot)$ defined in
\eqref{fn}. They coincide if $d=1$.
Clearly, for $d\geq 1$
$$
{\sf K}(r)\leq {\sf R}_X(r)\leq {\sf K}(r\sqrt{d}),\;\;\;\;r\in \mathbb{N}.
$$
Consequently, if ${\sf K}\in \mathcal{L}(\mathbb{N})$, then ${\sf R}_X\in \mathcal{L}(\mathbb{N})$,
and vise versa if
${\sf R}_X\in \mathcal{L}(\mathbb{N})$, then ${\sf K}\in \mathcal{L}(\mathbb{N})$.
Now for a sequence
$(C_r)_{r\in\mathbb{N}}$ it is not difficult to obtain the desired result following the scheme of the proof of
Theorem  \ref{th1} and
using ${\sf R}_X$ instead of ${\sf K}_X$. $\square$
\vskip0,2cm
\begin{Rm}\label{rm2} Lemma $\ref{lm2}$ shows that in Theorems $\ref{th1}$ and $\ref{th2}$
instead of normalizations $\sqrt{{\sf var} S_n}$ and
$\sqrt{{\sf var} S(C_r)}$ for partial sums one can use $\sqrt{\langle n \rangle \, K_X(n)}$
and $r^{d/2}\sqrt{K(r)}$, respectively.
\end{Rm}
\vskip1cm
{\bf  Acknowledgments.} The author is grateful to Professors I.Kourkova and G.Pages for invitation
to LPMA of the University Pierre and Marie Curie, he would like also to thank  all the members of the
LPMA for hospitality.

\end{document}